\documentclass[11pt,a4paper,reqno]{article}
\usepackage{caption}
\usepackage[english]{babel}
\usepackage{setspace}
\usepackage{float}
\usepackage[letterpaper,top=2cm,bottom=2cm,left=3cm,right=3cm,marginparwidth=1.75cm]{geometry}
\usepackage{amssymb,amsthm} 
\usepackage{url}
\usepackage{adjustbox}
\usepackage{multirow}
\usepackage{amsmath,tabularx}
\usepackage{empheq,etoolbox}
\patchcmd{\subequations}
  {\theparentequation\alph{equation}}
  {\theparentequation.\alph{equation}}
  {}{}

\usepackage{algorithm,algpseudocode}
\usepackage{subcaption}
\usepackage{color}
\usepackage[sort,nocompress]{cite}
\usepackage{booktabs}

\newcommand{\xbf}{\mathbf{x}}

\newcommand{\pbf}{\mathbf{p}}

\newcommand{\abf}{\mathbf{a}}
\newcommand{\zerobf}{\mathbf{0}}

\numberwithin{equation}{section}
\usepackage{graphicx}
\usepackage[colorlinks=true, allcolors=blue]{hyperref}

\usepackage{authblk}

\title{On Optimal Control at the Onset of a New Viral Outbreak}
\author[1]{Alexandra Smirnova\thanks{asmirnova@gsu.edu. Supported by NSF award 2011622 (DMS Computational Mathematics)}}
\author[1]{Xiaojing Ye\thanks{xye@gsu.edu. Supported by NSF award 2152960 (DMS CDS\&E) and 2307466 (DMS Applied Mathematics).}}
\affil[1]{Department of Mathematics \& Statistics,  Georgia State University, Atlanta, USA}

\begin{document}
\maketitle

\begin{abstract}

We propose  a versatile model with a flexible choice of control for an early-pandemic outbreak prevention when vaccine/drug is not yet available.
At that stage, control is often limited to non-medical interventions like social distancing and other behavioral changes.
 For the SIR optimal control problem, we show that the running cost of control satisfying  mild, practically justified  conditions generates  an optimal  strategy, $u(t)$, $t \in [0, T]$, that  is sustainable up until some moment $\tau \in [0 ,T)$. However,  for any $t \in [\tau, T]$,  the function  $u(t)$ will  decline as $t$ approaches $T$, which may cause the number of newly infected people to increase.  So, the window from $0$ to $\tau$ is the time for public health officials  to prepare alternative mitigation measures, such as vaccines, testing,  antiviral medications, and others. In addition to theoretical study, we develop a fast and stable computational method for solving the proposed optimal control problem.  The efficiency of the new method is illustrated with numerical examples of optimal control trajectories for various cost functions and weights. Simulation results provide a comprehensive demonstration of the effects of control on the epidemic spread and mitigation expenses, which can serve as invaluable references for public health officials.


\end{abstract}


\section{Introduction}

The circulation of infectious diseases, such as COVID-19,   is shaped by multiple parameters including control interventions \cite{P20},  environmental factors \cite{WM04},   immunity patterns \cite{OD21},  superspreading events \cite{L05},  and behavior changes \cite{10.1001/jamanetworkopen.2021.15959}.  These factors impact the early growth dynamics \cite{S04} and the basic reproduction number \cite{WE21}, which quantifies the number of secondary cases per primary case in a completely susceptible population.

Since the first COVID-19 case was detected in December 2019, the disease spread rapidly causing a worldwide pandemic. While some infected people experience only mild or moderate symptoms, others can get seriously ill and require immediate medical intervention \cite{TimeToDeatch, Shao, Wu, He}. Among high risk individuals are elderly people and those with underlying health conditions such as cancer, diabetes, chronic respiratory disease, and others \cite{CDCrisk}. As of April 14, 2024, there have been 775,335,916
 confirmed cases of COVID-19, including 7,045,569 deaths \cite{WHOpop}. Important factors contributing to the alarming rise in COVID‑19 cases at the early stage of the pandemic were high reproduction number, a large number of "silent spreaders" (especially among young people), and a relatively long incubation period \cite{CDCrisk1}. In the absence of vaccines and antiviral treatments in late 2019 and early 2020 \cite{HHStreatment}, mitigation measures  such as social distancing (including full or partial lockdowns),  restrictions on travel and mass gatherings, isolation and quarantine of confirmed cases, change from in-person to online education, and other similar tools emerged as the key ways of control and prevention \cite{M3, OECDfactors}.
While these measures proved to be effective in a short-term, they are hard to sustain  in a long run due to their negative impact on mental health coupled with high social and economic cost. Hence, since the start of COVID-19, balancing pros and cons of early non-medical interventions has come to the forefront (not only to contain COVID-19, but also to prepare for future epidemic outbreaks) \cite{L23, M1,M2,HADI2020102317,aronna2020model,Pazos2020.05.27.20115295,svoboda2022infection,yuan2022school,denoel2022decision,panovska2020determining,paltiel2020assessment,Barlow, S20, S21, S22a, LGC10}.

In this paper,  we consider an optimal control problem for  SIR compartmental model (Susceptible $\rightarrow$ Infectious $\rightarrow$ Removed)  of early disease transmission.   We design a running cost of control with mild, practically justified conditions that give rise to the optimal control strategy, $u(t)$, which does not exceed its admissible upper bound for the entire duration of the study period, $[0,T]$. Our theoretical analysis  indicates that at the early stage of an outbreak, the optimal control strategy, $u(t)$,   may be growing until some moment $\tau \in [0 ,T)$. However,  for any $t \in [\tau, T]$,  the function  $u(t)$ will  decline as $t$ approaches $T$, which may cause the number of newly infected people to increase.  So, the window from $0$ to $\tau$ is the time for public health officials  to prepare alternative mitigation measures, such as vaccines, testing,  antiviral medications, and others. Our theoretical findings are illustrated with important numerical examples showing optimal control trajectories for various cost parameters. To learn the optimal control function $u(t)$, we employed a deep learning based numerical algorithm, where $u(t)$ is parameterized as a deep neural network (DNN). The implementation, training and testing of all methods were conducted in Python 3.9.6 with PyTorch 2.1.0 and Torchdiffeq 0.2.3.

\section{A Strategy for Early Intervention}

At the onset of an emerging epidemic, in the absence of a vaccine and antivirals \cite{HHStreatment, S15, S18, S22}, the transmission of individuals between different stages of infection is often described by a classical SIR (Susceptible $\rightarrow$ Infectious $\rightarrow$ Removed) compartmental model \cite{kermack1927contribution, kudryashov2021analytical, S24}. For this early and relatively short phase, it is reasonable to infer that natural birth and death balance one another and, therefore, can be omitted. With the disease death rate varying between age and risk groups and being hard to estimate early on, the removed class is assumed to combine recovered and deceased people. Finally, due to the fast dynamic of the initial pre-vaccination stage, we suppose that recovered individuals develop at least a short-term immunity and don't move back to the susceptible class until the end of the study period. Under these assumptions, the SIR (Susceptible $\rightarrow$ Infectious $\rightarrow$ Removed/{\it Immune $+$ Deceased}) model is given by
 the following system of ordinary differential equations:
\begin{align}
\dfrac{d \mathcal{S}}{d t} &= -\beta \dfrac{\mathcal{S}(t) \mathcal{I}(t)}{N} \notag\\
\dfrac{d \mathcal{I}}{d t} &= \beta \dfrac{\mathcal{S}(t) \mathcal{I}(t)}{N}- \gamma \mathcal{I}(t)  \\
\dfrac{d \mathcal{R}}{d t} &= \gamma \mathcal{I}(t)    \notag
\end{align}
The primary goal of our study is to look at possible control strategies that can be effectively introduced at the early ascending stage of an outbreak before more robust mitigation measures, such as vaccines and viral medications, become available. The most common early mitigation measures, which were broadly used during the recent COVID-19 pandemic, include physical distancing, enhanced personal hygiene, mask wearing, awareness, and others. Their primary goal is to "flatten the curve", that is, to reduce the daily number of new infections and, as the result, to reduce the number of virus-related deaths.
The SIR model with enforced control, $u=u(t)$, and normalized dependent variables, $S(t):= \dfrac{ \mathcal{S}(t)}{ N}$, $I(t):= \dfrac{ \mathcal{I}(t)}{ N}$, and $R(t):= \dfrac{ \mathcal{R}(t)}{ N}$, takes the form $\frac{d{\bf x}}{dt} = f({\bf x},u)$, where
\begin{align} \label{2}
f_1({\bf x},u) &:= -\beta (1-u(t)) S(t) I(t)  \notag\\
f_2({\bf x},u) &:= \beta (1-u(t)) S(t) I(t)- \gamma I(t)  \\
f_3({\bf x},u) &:= \gamma I(t),    \notag
\end{align}
and ${\bf x} = [S, I, R]^\top$. In the above, the admissible set for the control function, $u=u(t)$, is assumed to be
\begin{align*}
    \mathcal{U} = \left\{ u \in \mathcal{L}^1[0,T], \quad 0 \le u(t) < 1 \right\}.
\end{align*}
In Lemma 2.1 below we show that following the introduction of a time-dependent transmission rate, $\beta(t):=\beta (1-u(t))$, the model $\frac{d{\bf x}}{dt} = f({\bf x},u)$ remains correct in the sense that the trajectories $(S(t),I(t),R(t))$ starting in a positive octant do not leave the octant and are defined for all $t>0$.

{\bf Lemma 2.1}.   {\it Let $u(t)$ be an admissible control trajectory with ${\bf x}(t)$ satisfying
$\frac{d{\bf x}}{dt} = f({\bf x},u)$ defined in (\ref{2}) and $$(S(0),I(0),R(0)) \in \Delta^{2}:=\{(z_1,z_2,z_3)\in \mathbb{R}^{3}: z_1+z_2+z_3=1,\ z_1,z_2,z_3\ge 0\},$$ the probability simplex in $\mathbb{R}^{3}$. Then  $(S(t),I(t),R(t)) \in \Delta^{2}$ for all $t\ge 0$.}

{\bf Proof.} We first notice that the solution to the system (\ref{2}) satisfies
\begin{align}
    S(t) & = S(0)e^{-\int_{0}^{t} \beta(1-u(s))I(s)ds} \label{eq:St}\\
    I(t) & = I(0)e^{\int_{0}^{t} (\beta(1-u(s))S(s)-\gamma)ds} \label{eq:It}\\
    R(t) & = R(0) +\gamma \textstyle\int_{0}^{t} I(s)ds \label{eq:Rt}
\end{align}
Therefore, $S(t),I(t) \ge 0 $ for all $t\ge 0$ due to \eqref{eq:St} and \eqref{eq:It}, respectively, where the latter further implies $R(t) \ge 0$ due to \eqref{eq:Rt}. Moreover, since $(S(t)+I(t)+R(t))' = 0$ for all $t$ due to the dynamics (\ref{2}), we know $S(t)+I(t)+R(t)=1$ for all $t\ge 0$. Combining these facts, we conclude that $(S(t),I(t),R(t)) \in \Delta^{2}$ for all $t\ge 0$. \hskip 70 mm $\Box$

One can easily see that model  (\ref{2}) yields
\begin{align} \label{2a}
\dfrac{d I}{d t} &= \left(\frac{\beta}{\gamma} \,(1-u(t)) S(t) - 1\right)\gamma I(t),
\end{align}
where $\beta/\gamma$ is the basic reproduction number. Clearly,  if  $\beta S(0)/\gamma<1$, then the virus is contained (even though it can still benefit from mitigation  measures that would further reduce the daily number of new infections).
It also implies that an obvious way of controlling the disease, should  $\beta S(0)/\gamma$ be greater than $1$, is to choose $u(t)$ such that  $\beta S(t)(1-u(t))/\gamma \ll 1$. That is, $1>u(t)\gg 1-\frac{\gamma}{S(t)\beta}$. However, all things considered, if the basic reproduction number, $\beta/\gamma$, is large, this kind of control may not be feasible. Indeed,  while the right interventions at the onset of the disease save lives and protect the health of the population, they come with social, psychological,  and economic costs. Therefore, policymakers have a difficult task of balancing the benefits to public health and the negative outcomes of their preventive measures. Mathematically, this comes down to solving the optimal control problem, where
the main goal  is to reduce the daily number of new infections, $\beta (1-u(t)) S(t) I(t)$, while also minimizing the cost of preventive measures, $\lambda \,c(u(t))$. This gives rise to the following objective functional:
\begin{align*}
    J({\bf x}, u): = \int_{0}^{T} \biggl\{(\beta(1-u(t)) S(t) I(t) + \lambda \,c(u(t))\biggr\}\,dt,\quad \lambda >0.
\end{align*}
According to (\ref{2}), this $J({\bf x}, u)$ can be written as
\begin{align} \label{3}
       J({\bf x}, u)  = S(0)- S(T) + \lambda\,\int_{0}^{T}   c(u(t))\,dt,\quad \lambda >0.
\end{align}
Evidently, the choice of the cost, $c(u(t))$,  has a major impact on the resulting control strategy. It is important to define $c=c(u)$ and the corresponding Lagrangian, denoted below as $L(x,u;p,q)$, in such a way that the optimal solution, $u=u(t)$, is guaranteed to take values between $0$ and $1$. In other words, it should never be a feasible strategy  for $u=u(t)$ to become negative, and the cost of control, $c(u(t))$, should get extremely high as $u(t)$ approaches $1$ (unless the regularization parameter, $\lambda$, is very small).

For various epidemic models, a very common choice of $c(u(t))$ is $c(u(t))=u^2(t)$ \cite{M1, M2}, which is often implemented in conjunction with the forward–backward sweep numerical algorithm for the computation of optimal control  $u=u(t)$ \cite{L07}.
However, as pointed out in \cite{M1, M2}, there are some drawbacks of setting $c(u(t))$ to $ u^2(t)$. Indeed, since this function has a finite penalty at $u(t)=1$, an explicit constraint $u(t) \le 1$ must be enforced. Without this constraint, it is easy to get $u(t)>1$ (especially for small values of $\lambda$), which results in unrealistic  strategy leading to  $S'(t)>0$. And even with the constraint $u(t) \le 1$, the optimal control function often reaches the nonviable "ultimate" level, $u(t)=1$, for the better part of the study period.

To avoid the above scenario, in our theoretical and numerical analysis, we consider $c=c(u)$ such that  $\lim_{u\to 1^-}c(u)=\infty$. This important requirement, along with the assumption that the cost, $c=c(u)$, is twice continuously differentiable in its domain containing $[0,1)$, with $\frac{d^2 c}{d u^2}>0$, $c(0)=0$, $c'(u)> 0$ for $u\ge 0$, $c'(u)<0$ for $u<0$, allows us to design an optimal control problem leading to  $u(t)<1$ during the entire study period, $[0,T]$. For numerical simulations,  we employ and compare $4$ different cost functions for the optimal control strategy
$$
c_1(u) = -  \ln(1-u^2), \,\,
    c_2(u) = -  u\ln(1-u), \,\,
    c_3(u)  = -(u + \ln(1-u)), \,\,\mbox{and}\,\,
    c_4(u)  =  u^2.
$$
All these functions, except for $c_4(u)$, have infinite penalty at $u(t)=1$, and the function $c_4(u)  =  u^2$ is used for comparison.

\section{Properties of Optimal Control}

In what follows, we prove our main theoretical result, Theorem 3.1,  that has major practical implications. Namely, we show that beginning with some moment, $\tau \in [0,T)$, the optimal control strategy, $u=u(t)$, introduced in the previous section, is declining, which may cause the number of newly infected people to increase.  So, the window from $0$ to $\tau$ is the time for public health officials  to prepare alternative mitigation measures, such as vaccines, testing,  antiviral medications, and others.

{\bf Theorem 3.1}. {\it Assume that $u=u(t)$ is an optimal control trajectory for the objective functional (\ref{3}) constrained by  the system
$\frac{d{\bf x}}{dt} = f({\bf x},u)$ defined in (\ref{2}) and  by the inequality  $u(t)\ge 0$ for all $t\in [0,T]$. Let $c=c(u)$ be twice continuously differentiable in its domain containing $[0,1)$, with $\frac{d^2 c}{d u^2}>0$, $c(0)=0$, $c'(u)> 0$ for $u\ge 0$, $c'(u)<0$ for $u<0$, and $\lim_{u\to 1^-}c(u)=\infty$. Then
there is  $\tau \in [0 ,T)$ such that for any $t \in [\tau, T]$,  the derivative, $u'(t)$, exists and $u'(t) \le 0$.}

{\bf Proof.}  Without loss of generality, we assume $S(t),I(t)>0$ (since otherwise it is clear that $u'(t)=0$). Given the constraints $\frac{d{\bf x}}{dt} = f({\bf x},u)$ and $u(t)\ge 0$ for all $t\in [0,T]$, the optimal control problem  results in the following Lagrangian
\begin{align}
L({\bf x},u;p,q)& = S(0)-S(T) + \int_0^T \left\{\lambda c(u(t)) -  p(t)\cdot ({\bf x}'(t) - f({\bf x}(t),u(t)))  -  q(t) u(t)\right\} dt \notag\\
&- p(0)\cdot ({\bf x}(0)-{\bf x}_0), \quad p(t) = [p_1(t), p_2(t), p_3(t)]^\top.
\end{align}
Then the Karush–Kuhn–Tucker (KKT)  conditions are as follows
\begin{align}
   (C1)\quad & \lambda c'(u) + p \cdot \partial_u f({\bf x},u) - q = 0 \\
   (C2)\quad & p' = -\partial_x f({\bf x},u)^{\top} p, \quad p(T) = [-1,0,0]^\top \\
   (C3)\quad & {\bf x}' = f({\bf x},u), \quad {\bf x}(0) = {\bf x}_0 \\
   (C4)\quad & q(t) \ge 0, \quad u(t) \ge 0, \quad q(t)u(t) = 0, \quad \forall t\in [0,T].
\end{align}
From (C1) we conclude  that $\lambda c'(u) - q = - p \cdot \partial_u f({\bf x},u) = -\beta SI(p_1-p_2),$ which is differentiable and hence continuous since all terms on the right-hand side are so due to (C2) and (C3). Furthermore, since $p_1(T)=-1<0=p_2(T)$ and $S(t),I(t)>0$,  there exist $\tau_1,\tau_2 \in [0,T)$ such that $p_1(t) < p_2(t)$ for all $t \in [\tau_1,T)$ and $p_1(t) < 0$ for all $t \in [\tau_2,T)$. Let $\tau = \max(\tau_1,\tau_2)$, then
\begin{equation}\label{p1a}
p_1(t) < 0\quad\mbox{and}\quad \lambda c'(u(t))-q(t)= -\beta S(t)I(t)(p_1(t)-p_2(t)) >0 \quad \forall t\in [\tau, T].
\end{equation}
Now we restrict our discussion to $[\tau,T]$. If $c'(u(t))=0$ for some $t$, then $u(t)=0$ by the property of $c$, and hence $q(t) = \beta S(t)I(t)(p_1(t)-p_2(t)) < 0$ which contradicts to (C4). Therefore $c'(u(t)) > 0$ for all $t$ by the assumptions on $c$ and $0<u(t)<1$, where the latter also implies $q(t)=0$ for all $t$ due to (C4). In summary, over $[\tau, T]$, we have $$\lambda c'(u(t)) + \beta S(t)I(t)(p_1(t)-p_2(t))=0. $$  Then by implicit function theorem we know $u'$ exists and
\begin{align}\label{est1}
u'(t) = -\frac{\beta [S(t)I(t)(p_1(t)-p_2(t))]'}{\lambda c''(u(t))} \quad \mbox{for all}\quad  t \in [\tau, T].
\end{align}
Taking into consideration (\ref{2}) and (C2), for all $t \in [\tau, T]$, one has $\,p_3(t) = 0\,$ and
\begin{equation}\label{p}
    \begin{cases}
         \dot{p_1} = -\beta\,(p_2 - p_1)(1-u)  I \\
         \dot{p_2} = -\beta\,(p_2 - p_1)(1-u) S + \gamma \,p_2  \\
         p_1(T) = -1, \,\, p_2 (T) = 0.
    \end{cases}
\end{equation}
From system (\ref{p}), one concludes
\begin{equation}\label{p_dif}
 \dot{p_2}- \dot{p_1}  = \beta\,(p_2 - p_1)(1-u) (I-S)+ \gamma \,p_2.
\end{equation}
Combining (\ref{2}) and (\ref{p_dif}),  one can rewrite $[S(t)I(t)(p_1(t)-p_2(t))]'$ as follows
\begin{align*}
[S(t)I(t)(p_1(t)-p_2(t))]'&=\beta\,(p_2 - p_1)(1-u) (I-S)SI+\gamma p_2SI\\
&+(p_2 - p_1)\bigl\{-\beta(1-u)SI^2+\beta(1-u)S^2I-\gamma SI)\bigr\}\\
&=\bigl\{\beta\,(p_2 - p_1)(1-u) (I-S)+\gamma p_2+\beta\,(p_2 - p_1)\\
&\times(1-u) (S-I)-\gamma(p_2-p_1)\bigr\}SI=\gamma p_1 SI.
\end{align*}
According to (\ref{p1a}) and (\ref{est1}), this yields for all $t \in [\tau, T]$,
\begin{align}\label{est}
u'(t) = -\frac{\beta [S(t)I(t)(p_1(t)-p_2(t))]'}{\lambda c''(u(t))} = \frac{\beta\gamma S(t)I(t)p_1(t)}{\lambda c''(u(t))} < 0.
\end{align}
This completes the proof. \hskip 110 mm $\Box$\vskip 2mm

{\bf Remark 3.2}  As it follows from (\ref{est}), the impact of $u(t)$ scaling down  towards the end of the early stage of the outbreak will depend on the weight, $\lambda$. If the wight on control is relatively high (see Figures 3 and 4 in Section 5), then the decline in $u(t)$ for $t\ge \tau$ can be substantial, which will result in a notable surge in the daily number of infected individuals, $\mathcal{I}(t)$, for $t\ge \tau$. On the other hand, if the weight, $\lambda$,  is small (as in Figure 5), then  by the time $t=\tau$ the epidemic is effectively contained. Hence  the decline in $u(t)$ for $t\ge \tau$ is negligible and the daily number of infected people, $\mathcal{I}(t)$, remains very low for $t\ge \tau$.

 Note that cost functional (\ref{3}) aims to minimize the cumulative number of cases during the early phase of the disease, i.e., for $t\in [0,T]$. However, it does not guarantee that on any given day, the number of infected individuals in the optimally controlled environment is less than the number of infected individuals in the same environment but with  no control. As our experiments below illustrate, when  the weight on control, $\lambda$, is relatively high, towards the end of the study period in a controlled environment the daily number of infected humans, $\mathcal{I}(t)$, can potentially bypass the corresponding $\mathcal{I}(t)$ in the environment with no control (see Figures 3 and 4 in Section 5).

 The objective functional (\ref{3}) is set to minimize the cumulative number of infections, $S(0)-S(T)$, while keeping the negative impact of mitigation measures at bay. This is achieved, for the most part, by reducing the daily number of new infections but also, apparently, by delaying some infections. On the bright side, $\mathcal{I}(t)$ in the controlled environments gets bigger than $\mathcal{I}(t)$ in the uncontrolled case only when $t$ is approaching $T$. It is reasonable to assume that at this time additional intervention measures become available that will gradually replace the initial set of controls.

In the running cost, rather than minimizing the daily number of new infections, one can also minimize the daily number of infected individuals. This gives rise to the following objective functional
\begin{align}\label{JI}
    \tilde J({\bf x}, u): = \int_{0}^{T} \bigl\{I(t) + \lambda \,c(u(t))\bigr\}\,dt= \frac{R(T)-R(0)}{\gamma} + \lambda \,\int_{0}^{T}c(u(t))\,dt.
\end{align}
That is, instead of maximizing $S(T)$, this functional aims to minimize $R(T)$.
 Using the similar argument as above, one  can show that this optimal control strategy, $\tilde u(t)$, will also be decreasing starting with some point $\tilde \tau\in [0,T)$.

\section{Numerical Algorithm for Learning Control}

To learn the optimal control function, $u:\,[0,T]\to \mathbb{R}$, which is guaranteed to take values in $[0,1)$ according to Theorem 3.1 above, we employ a deep learning based numerical algorithm. This algorithm can be easily modified to the case of vector-valued controls. At the first step, we parameterize $u$ as a deep neural network (DNN), denoted by $u_\theta$, with parameters $\theta \in \mathbb{R}^m$. In our experiments, we chose a simple fully connected network with both input and output layer dimension 1 (because  in our setting, the input is time $t \in[0,T]$ and the output is a scalar). We set $u_\theta$ to have 4 hidden layers and each layer is of size $10$.
Specifically, the function $u_{\theta}$ is defined by
\[
u_{\theta}(t) = w_{5}^{\top} \sigma(W_{4}z_{4} + b_4) + b_5, \quad \mbox{where} \quad z_{l+1} = \sigma(W_{l}z_{l}+ b_{l})\quad \mbox{for}\quad l=0,1,2,3
\]
and $z_0 := t \in \mathbb{R}$. Here $\sigma(z):=\tanh(z)$ is the activation function that applies to the argument componentwisely, and $\theta=(w_5,W_4,W_3,W_2,W_1,W_0,b_5,b_4,b_3,b_2,b_1,b_0)$ is a column vector that contains all the components of these variables, where $w_5 \in \mathbb{R}^{10}$, $W_4,W_3,W_2,W_1 \in \mathbb{R}^{10\times 10}$, $W_0 \in \mathbb{R}^{10\times 1}$, $b_5\in \mathbb{R}$, $b_4,b_3,b_2,b_1,b_0\in \mathbb{R}^{10}$, and all vectors are considered as column vectors. The number $m$ is the dimension of $\theta$ (i.e., the total number of components in $\theta$), which is $m=471$ in our case.
Introduce the notation
$\ell(\theta):= J(\mathbf{x}, u_\theta),$
where $J$ is defined in \eqref{3} and $\xbf$ follows the dynamics \eqref{2} with the given initial state $\xbf(0)=(S(0),I(0),R(0))$. To find the optimal $u_\theta$, we essentially need to compute $\nabla_{\theta} \ell(\theta)$ for any $\theta$ and apply the gradient descent to update $\theta$. In our algorithm, we employ the neural ordinary differential equation (NODE) method \cite{chen2018neural} which computes $\nabla_\theta \ell(\theta)$ in the following way.
First, with the given $\theta$, one solves the ODE forward in time:
\begin{equation}
    \label{eq:node-forward}
    \dot{\xbf}(t) = f(\xbf(t),u_\theta(t)), \quad t \in [0,T],
\end{equation}
with initial value $\xbf(0)=\xbf_0$ and $f(\xbf,u)$ defined in (\ref{2}).
Second, one solves the augmented adjoint equation backward in time:
\begin{equation}
    \label{eq:node-backward}
    (\dot{\pbf}(t),\ \dot{\abf}(t)) = - (\pbf(t)\, \partial_{\xbf} f(\xbf(t),u_\theta(t)),\  \pbf(t)\, \partial_u f(\xbf(t),u_\theta(t))\, \partial_\theta u_\theta(t)), \quad t \in [0,T],
\end{equation}
with terminal value $(\pbf(T),\xbf(T))=(-\nabla h(\xbf(T)), \zerobf)$. Here $\xbf,\pbf,\abf$ are all row vectors at each time $t$. Then it can be shown that $\nabla_\theta \ell(\theta)=\abf(0)$ \cite{chen2018neural}.  The algorithm is summarized in Algorithm \ref{alg:node} below. The implementation, training and testing were conducted in Python 3.9.6 with PyTorch 2.1.0 and Torchdiffeq 0.2.3. We initialize the parameter $\theta$ using Xavier initialization built in the PyTorch package.
\vskip 2mm
\begin{algorithm}
\caption{Neural ODE method to solve the optimal control problem of $u$}
\label{alg:node}
\begin{algorithmic}
\Require{Cost function $c$ and weight $\lambda>0$. Network structure $u_\theta$. Initial guess $\theta$.}
\Ensure{Optimal control $u_\theta$ with trained $\theta$.}
\Repeat
    \State Solve $\xbf$ forward in time using \eqref{eq:node-forward}.
    \State Solve $(\pbf,\abf)$ backward in time using \eqref{eq:node-backward}.
    \State $\theta \leftarrow \theta - \mu \abf(0)$.
\Until{converged.}
\end{algorithmic}
\end{algorithm}

A few  details about the performance of Algorithm \ref{alg:node} and our numerical simulations:
\begin{itemize}
\item In our experiments,  we try $4$ different cost functions, $c=c(u)$. Details and discussion will be given in Section 5;
\item  The weight, $\lambda$, scales the cost function, $c=c(u)$, and can be critical to the optimal control solution. In Section 5, we conduct empirical study on different values of $\lambda$;
\item One can apply any numerical integrator (e.g., Euler, mid-point, Runge-Kutta) to solve \eqref{eq:node-forward} and \eqref{eq:node-backward}. We used the 4th order Runge-Kutta method (rk4) built in the PyTorch package;
 \item One can either solve $(\pbf,\abf)$ as in \eqref{eq:node-backward} or solve $\xbf$ jointly backward in time to avoid saving $\xbf(t)$ obtained forward in \eqref{eq:node-forward} in the memory;
  \item In Algorithm \ref{alg:node}, one can choose the step size $\mu>0$ and terminate the algorithm after a prescribed number of iterations, $K$. We used Adam \cite{kingma2015adam:} with deterministic gradient and set these parameters as $\mu=0.001$ and $K=1,000$ and other parameters as default. The results appear to be stable for values around them;
 \item Since the control problem is not convex in $(\xbf,u)$, it is not guaranteed that our solution is the global minimizer. This is, unfortunately, a common issue in solving optimal control problems. Nevertheless, in all experiments, the numerical solutions obtained by our method appear to satisfy the constraint $u(t)\in [0,1)$ for all $t\in [0,T]$  and $du(t)/dt < 0$ starting with some $\tau \in [0 ,T)$;
     \item In the above Algorithm, we followed the idea of neural ODE \cite{chen2018neural} and set $\mathbf{a}(t)$ to be the auxiliary variable in order to compute the gradient of the loss function $\ell(\theta)$ with respect to $\theta$. Specifically, by solving the augmented adjoint dynamics (\ref{eq:node-backward}) backward in time, one can show that $\mathbf{a}(0)=\nabla_{\theta} \ell(\theta)$, which allows  to find the (local) minimizer of the loss function $\ell(\theta)$ by the gradient descent method. More details about the derivation can be found in  \cite{chen2018neural}.
   \end{itemize}

\section{Numerical Results and Discussion}
\begin{figure}[ht]
    \centering
    \includegraphics[width=.47\textwidth]{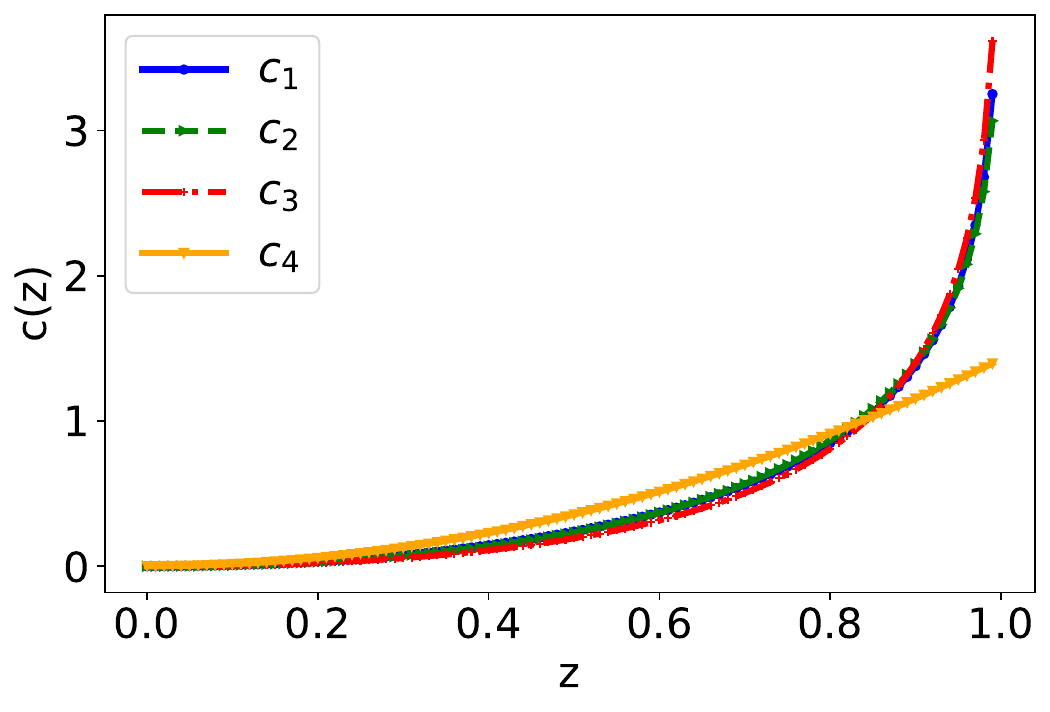}
    \caption{Comparison of different cost functions $c$.}
    \label{fig:compare_c}
\end{figure}

In this section, we apply the deep learning based numerical algorithm to solve the optimal control problem (\ref{3})
     subject to  SIR model (\ref{2}) with the following four cost functions:
\begin{align*}
    c_1(u) & = - 0.830071 \, \ln(1-u^2) \\
    c_2(u) & = - 0.672850 \, u\ln(1-u) \\
    c_3(u) & = - u - \ln(1-u) \\
    c_4(u) & = 1.424546 \,u^2.
\end{align*}
The weight $0.830071$ in $c_1$ is chosen to minimize the distance
\[
\int_0^1 w(z) |c_1(z) - c_3(z)|^2 \, dz,\quad w(z) = \sqrt{1-z^2},
\]
(the same for $c_2,c_4$). Doing so makes $c_i$'s close in the $w$-weighted 2-norm sense. See the comparison of these cost functions in Figure \ref{fig:compare_c}.

\begin{figure}[h]
    \centering
    \includegraphics[width=.44\textwidth]{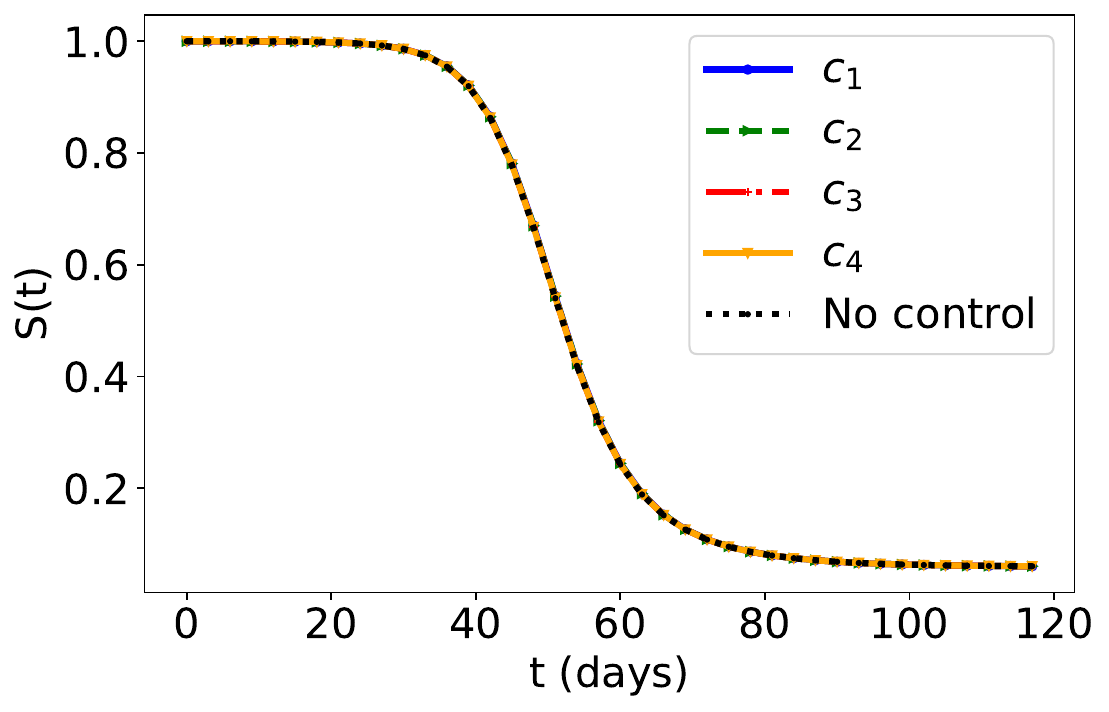}
    \includegraphics[width=.44\textwidth]{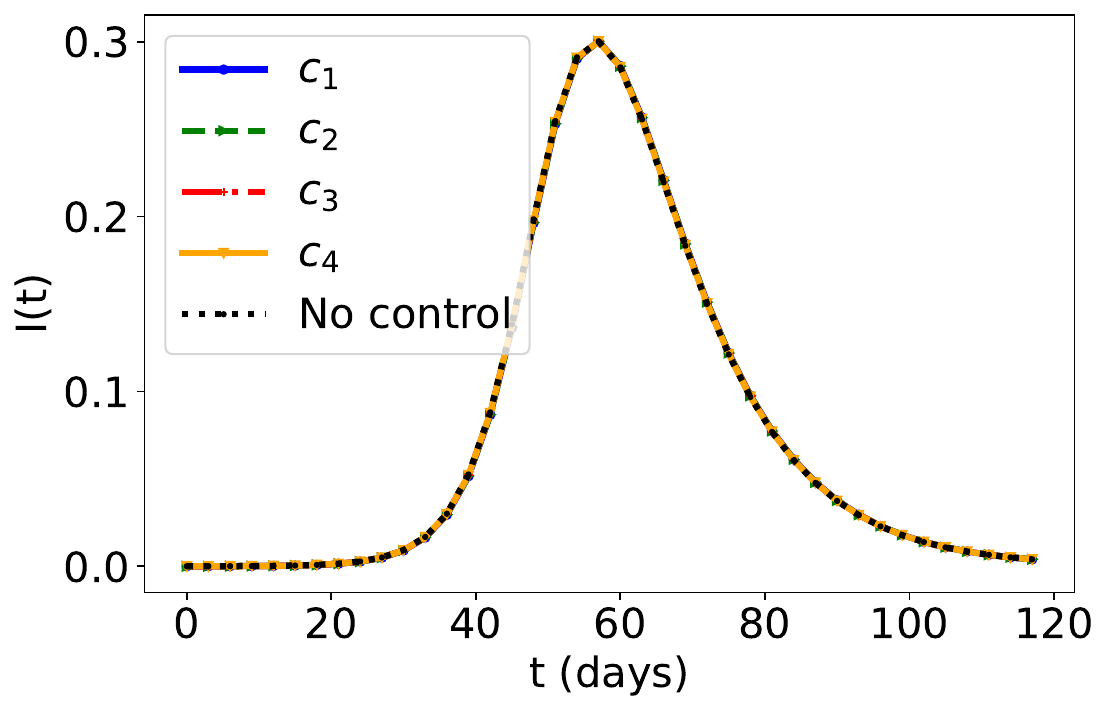}
    \includegraphics[width=.44\textwidth]{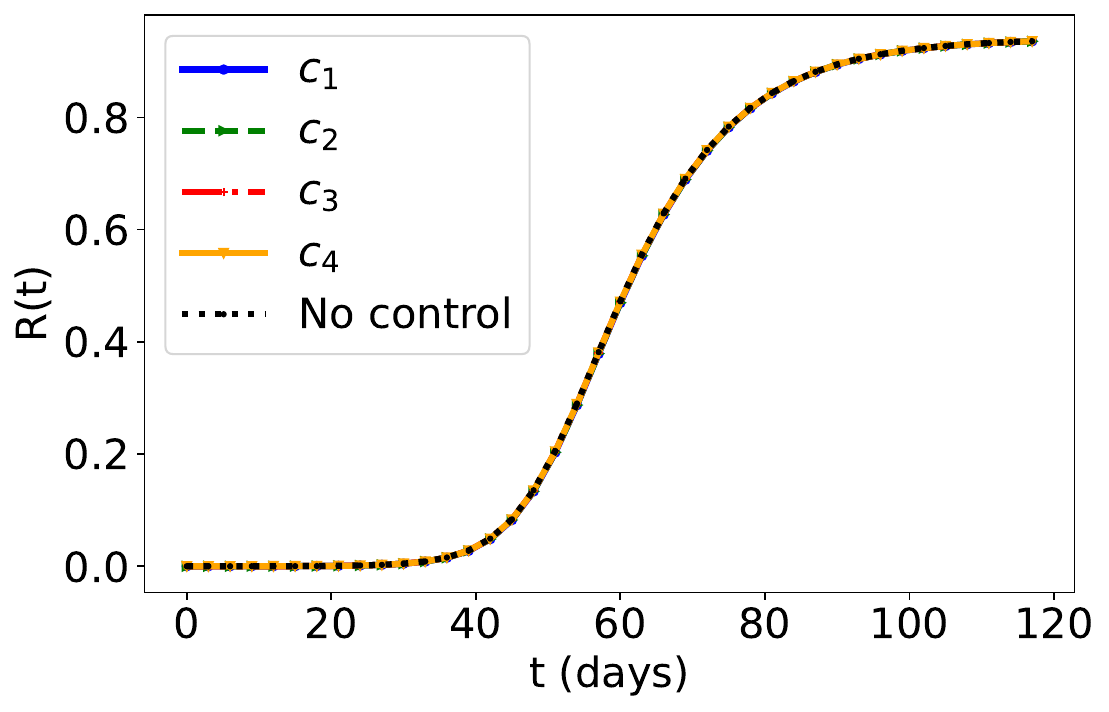}
    \includegraphics[width=.44\textwidth]{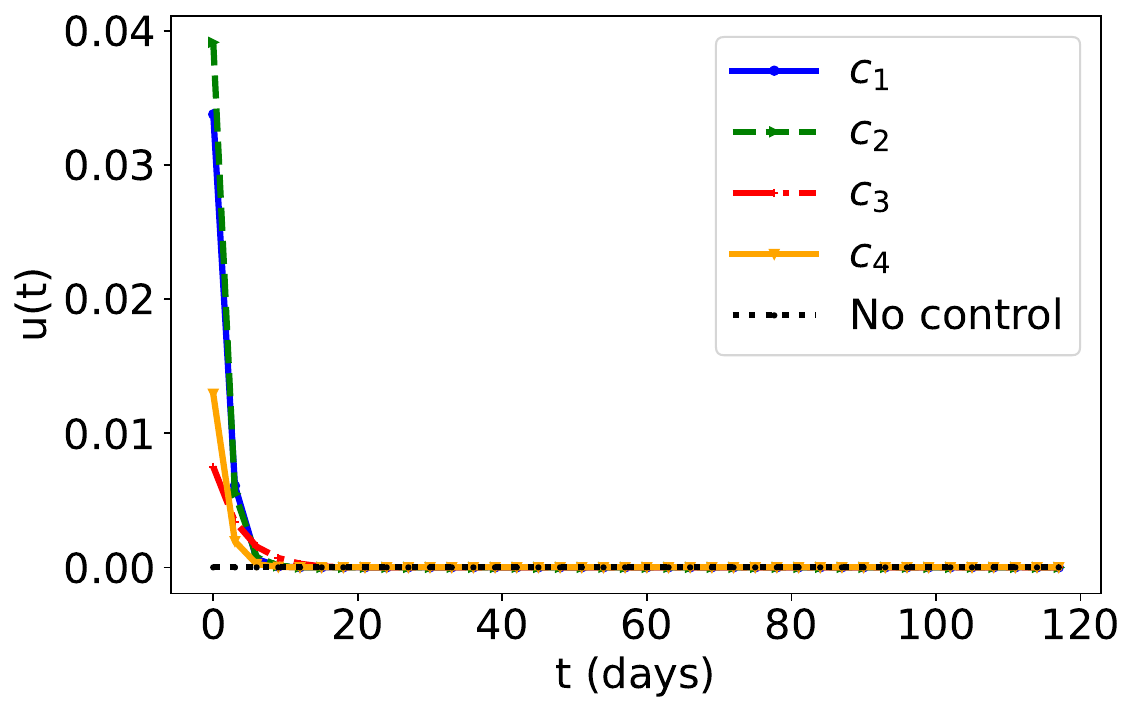}
    \caption{Weight $\lambda=0.1$.}
    \label{fig:lambda1e-1}
\end{figure}

In our experiments, we consider $\lambda$ as $0.1$, $0.05$, $0.01$, and $10^{-7}$. For $\lambda=0.1$, the environment is close to no control as shown in Figure \ref{fig:lambda1e-1}, since the penalty on control is weighted highly. All values of $\lambda$ larger than $0.1$ resulted in a similar behavior and hence they were omitted here.

\begin{figure}[h]
    \centering
    \includegraphics[width=.44\textwidth]{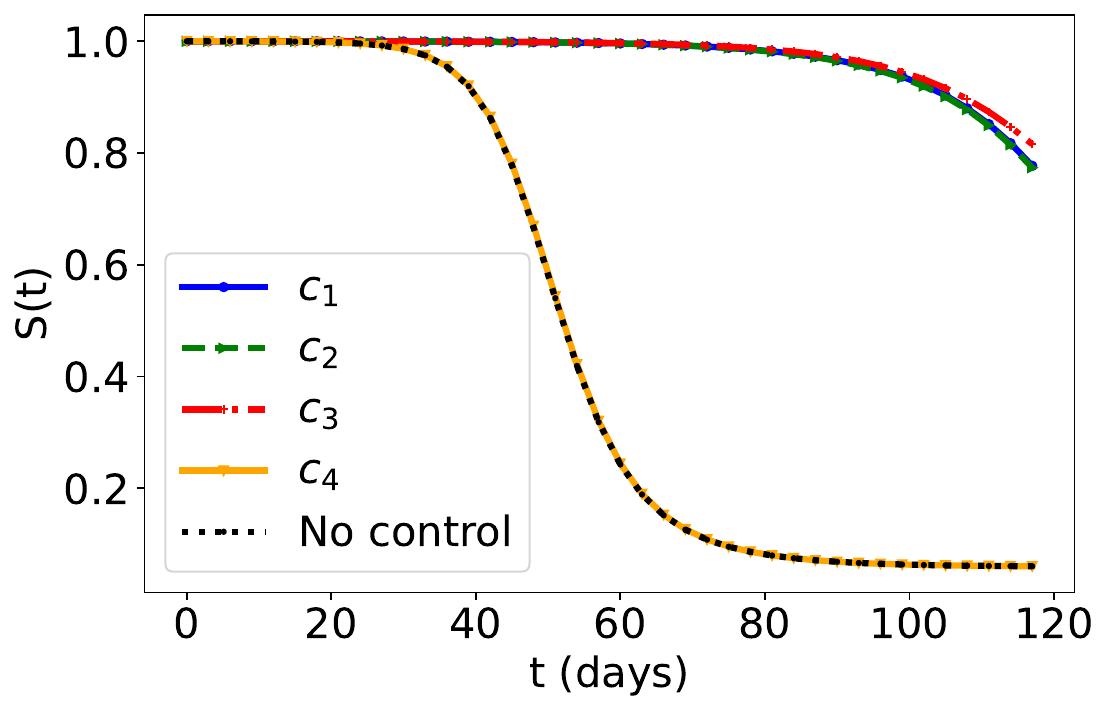}
    \includegraphics[width=.44\textwidth]{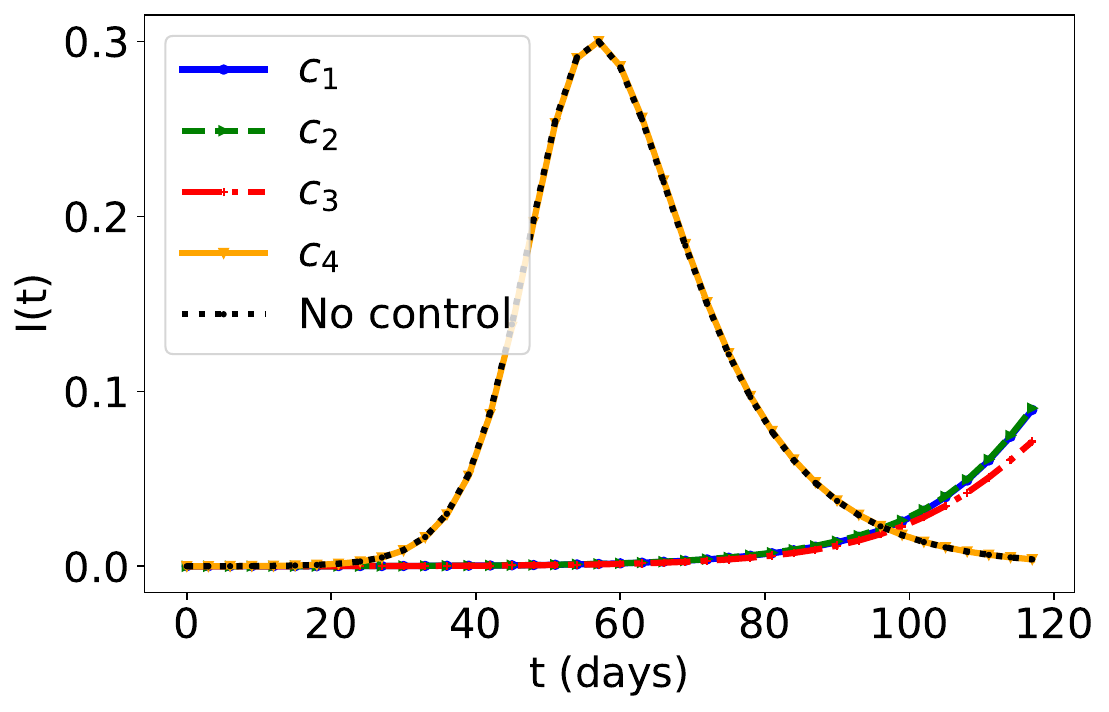}
    \includegraphics[width=.44\textwidth]{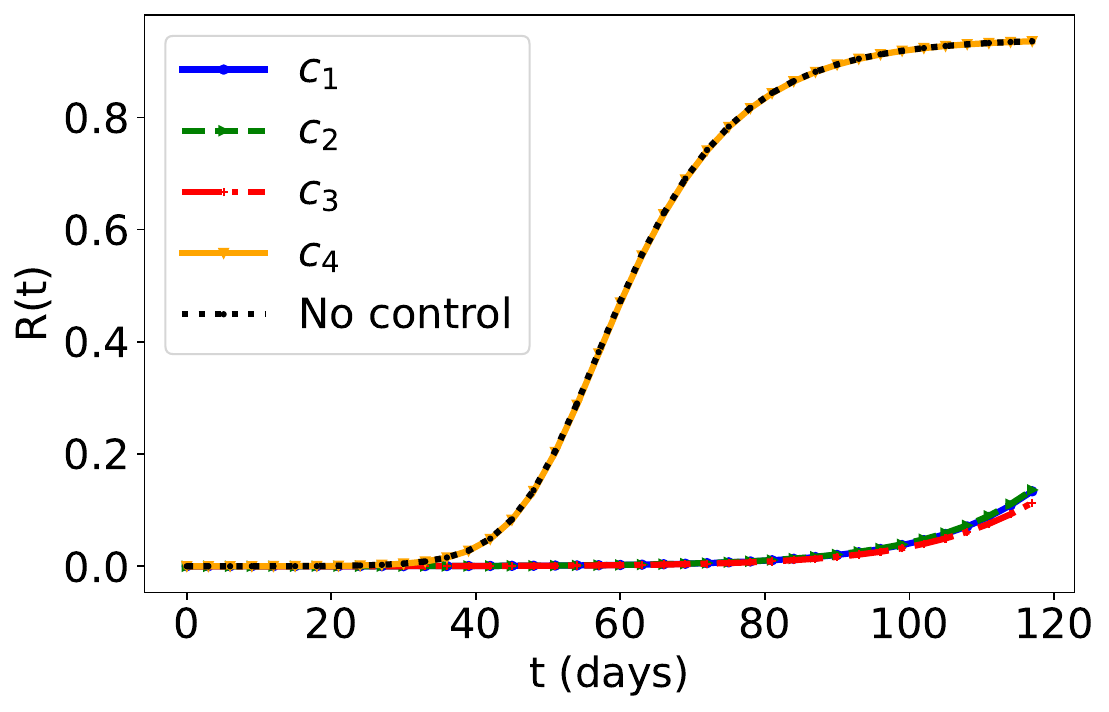}
    \includegraphics[width=.44\textwidth]{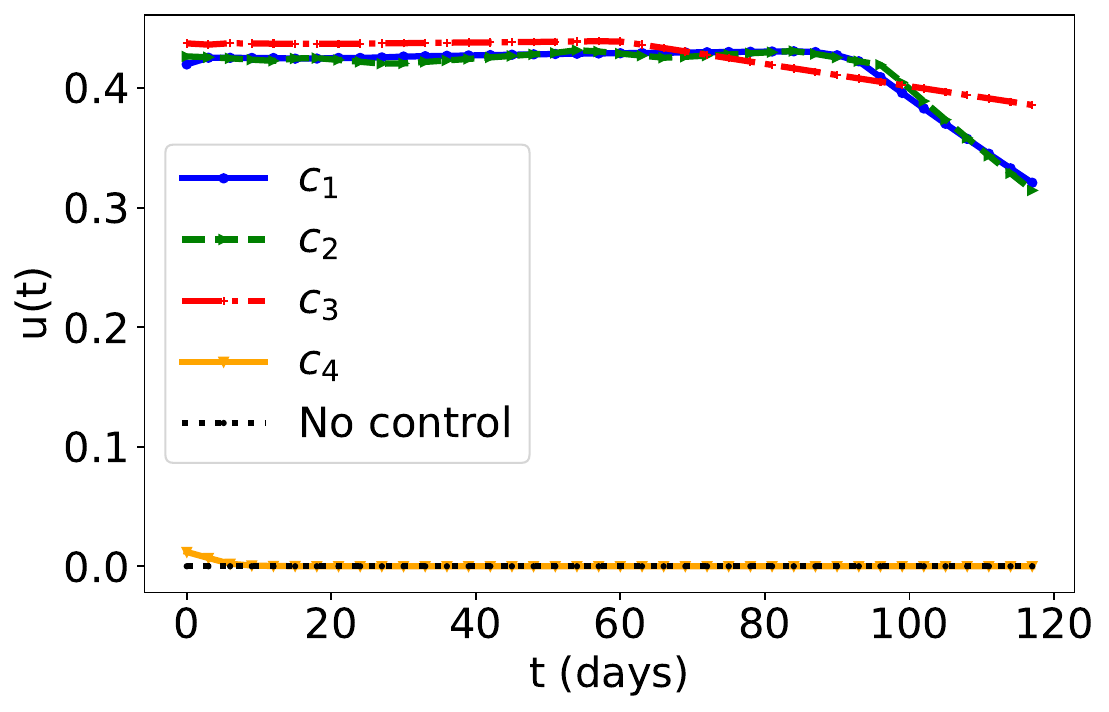}
    \caption{Weight $\lambda=0.05$.}
    \label{fig:lambda5e-2}
\end{figure}

As $\lambda$ gradually decreases, we see that controls start to make impact as they become cheaper to implement.  For example, as shown in Figure \ref{fig:lambda5e-2}, when $\lambda=0.05$, we observe that $c_1(u)$, $c_2(u)$, and $c_3(u)$ generate similar control strategies that effectively suppress the cumulative number of infected people (or equivalently maximize $S(T)$).

We note that all three controls, $c_1(u)$, $c_2(u)$, and $c_3(u)$, make a considerable positive impact on how the epidemic unfolds. Even though $\mathcal{I}(t)$ in the controlled environment is still growing (see Table 1), the daily number of infected people remains low  for a long time.  However, as expected from our theoretical analysis, for all three cost functions, $c_1(u)$, $c_2(u)$, and $c_3(u)$, and $\lambda =0.05$, the corresponding control strategies, $u(t)$,  begin to decrease after some point.

\begin{table}[h]
 \caption{Number of infected people $\mathcal{I}(t)$ versus day $t$ for $\lambda=0.05$}
    \centering
    \begin{tabular}{rrrrrr}
    \toprule
       Day & No control & $c_1$ & $c_2$ & $c_3$ & $c_4$ \\
       \midrule
0 & 200 & 200 & 200 & 200 & 200 \\
10 & 1626 & 423 & 423 & 407 & 1600 \\
20 & 12542 & 889 & 891 & 824 & 12341 \\
30 & 92164 & 1856 & 1882 & 1660 & 90714 \\
40 & 643531 & 3898 & 4004 & 3370 & 634826 \\
50 & 2358721 & 8104 & 8348 & 6785 & 2344806 \\
60 & 2852657 & 16674 & 17093 & 13546 & 2858895 \\
70 & 1722500 & 34405 & 35574 & 27592 & 1731405 \\
80 & 834814 & 69687 & 72449 & 57021 & 839919 \\
90 & 373748 & 137590 & 143188 & 118394 & 376160 \\
100 & 164980 & 276638 & 282921 & 245423 & 166065 \\
110 & 71628 & 567049 & 572008 & 479114 & 72103 \\
\bottomrule
    \end{tabular}
    \label{tab:lambda0.05}
\end{table}

Thus, towards the end of the study period in a controlled environment the daily number of infected humans, $\mathcal{I}(t)$, bypasses the corresponding $\mathcal{I}(t)$ in the environment with no control.

As mention in Remark 3.4, the objective functional (\ref{3}) is set to minimize the cumulative number of infections, $\mathcal{S}(0)-\mathcal{S}(T)$, which it does successfully as it is evident from Table 2.  This is achieved, for the most part, by reducing the daily number of new infections but also, apparently, by delaying some infections.

\begin{table}[h]
 \caption{Total number of infected people up to day $t$,  $N-\mathcal{S}(t)$, for $\lambda=0.05$}
    \centering
    \begin{tabular}{rrrrrr}
    \toprule
       Day & No control & $c_1$ & $c_2$ & $c_3$ & $c_4$ \\
       \midrule
0 & 200 & 200 & 200 & 200 & 200 \\
10 & 2339 & 734 & 731 & 709 & 2303 \\
20 & 18723 & 1835 & 1843 & 1724 & 18425 \\
30 & 138618 & 4154 & 4182 & 3786 & 136445 \\
40 & 993903 & 9034 & 9212 & 7991 & 980029 \\
50 & 4179153 & 19132 & 19626 & 16404 & 4146211 \\
60 & 7574716 & 39793 & 40813 & 33135 & 7557842 \\
70 & 8800810 & 82846 & 85337 & 67540 & 8795738 \\
80 & 9181982 & 169453 & 175482 & 138477 & 9180310 \\
90 & 9316270 & 339228 & 352155 & 285376 & 9315647 \\
100 & 9367572 & 682997 & 703593 & 592907 & 9367317 \\
110 & 9388987 & 1389004 & 1413628 & 1187340 & 9388880 \\
\bottomrule
    \end{tabular}
    \label{tab:lambda0.05_N-S}
\end{table}

As one can clearly see from Figure 3 and Tables 1 and 2, the cost of the optimal control strategy, $u(t)$, corresponding to $c_4$, is still very high for $\lambda =0.05$, and this control does not defeat  the outbreak. The reason for this control being different from $c_1(u)$, $c_2(u)$, and $c_3(u)$ can be understood from Figure 1, which shows that the cost, $c_4(u)$, is greater than the cost of all other controls between $u=0.2$ and $u=0.8$.

\begin{figure}[h]
    \centering
    \includegraphics[width=.44\textwidth]{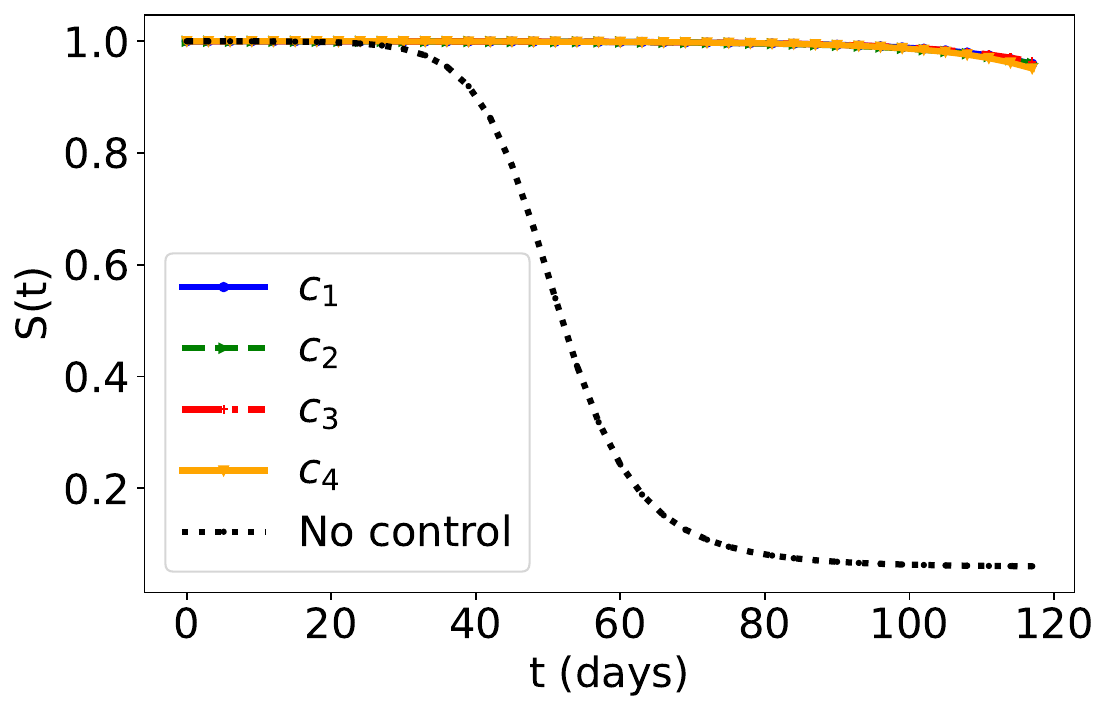}
    \includegraphics[width=.44\textwidth]{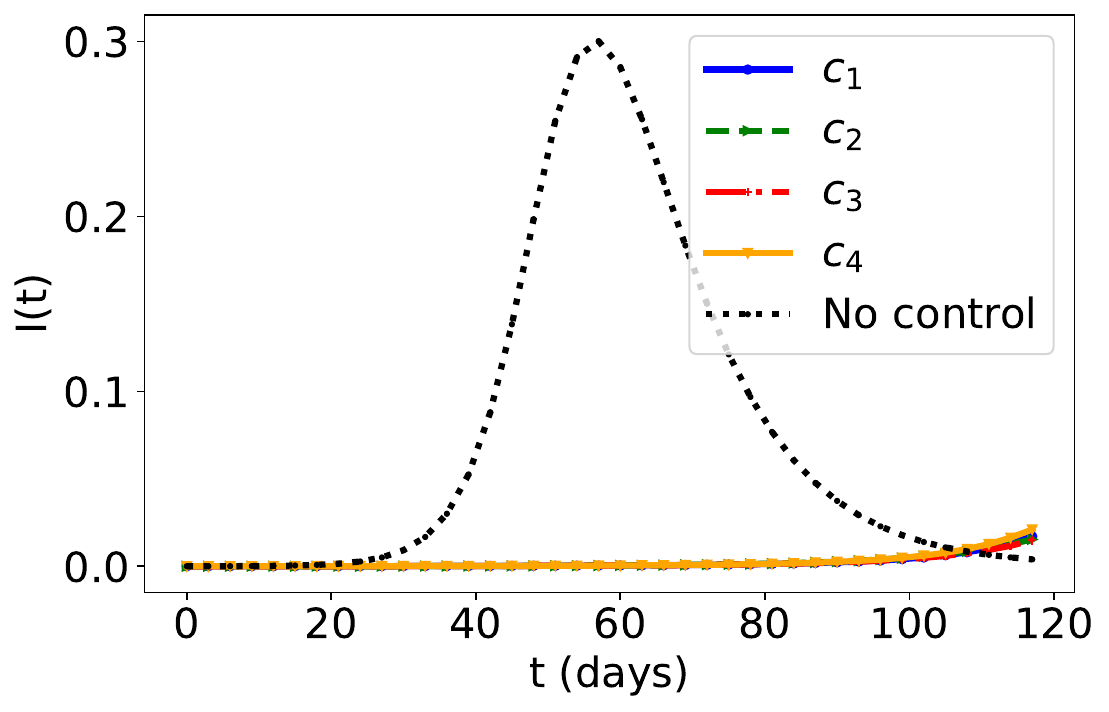}
    \includegraphics[width=.44\textwidth]{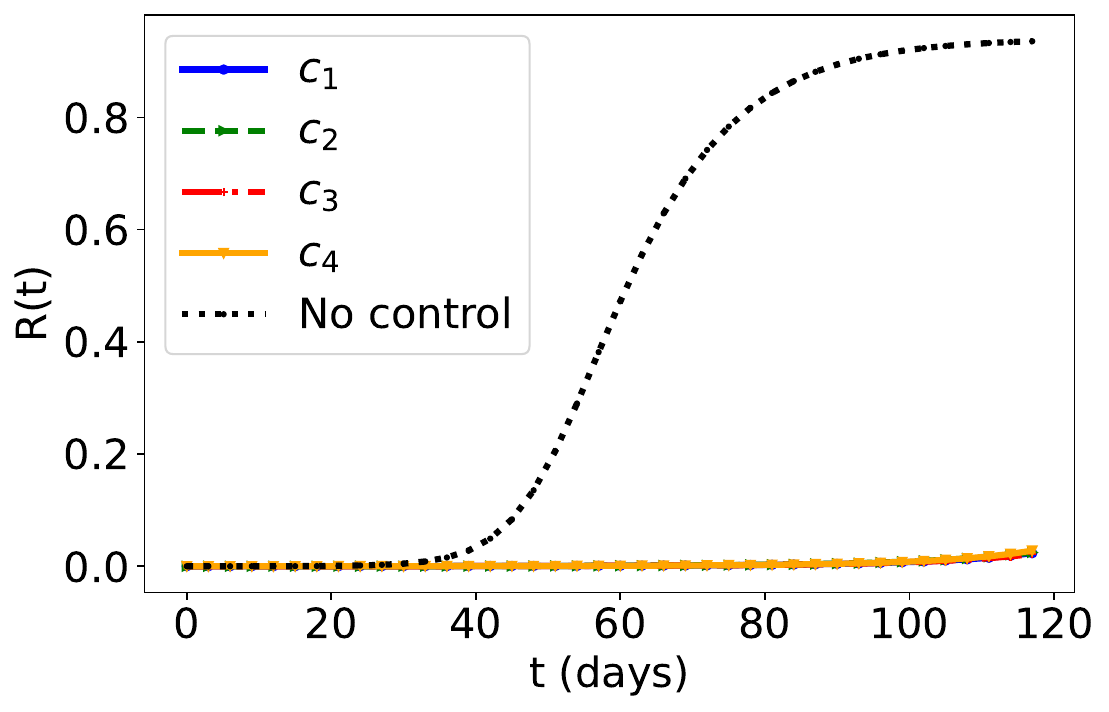}
    \includegraphics[width=.44\textwidth]{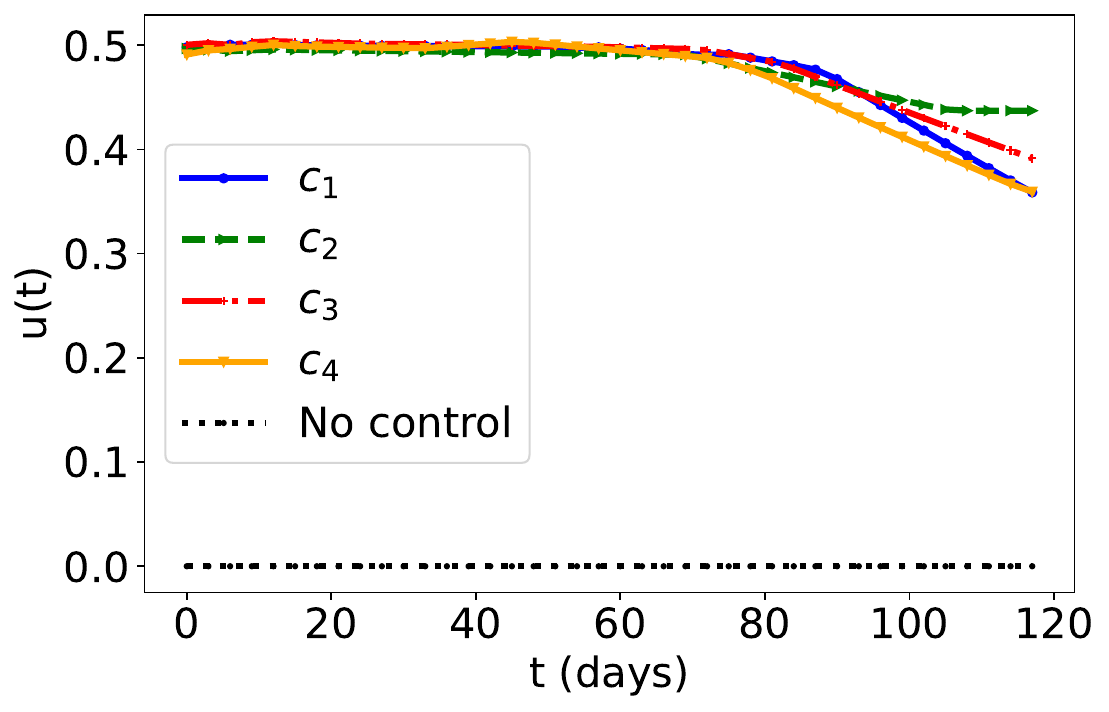}
    \caption{Weight $\lambda=0.01$.}
    \label{fig:lambda1e-2}
\end{figure}

When $\lambda$ reaches $0.01$, the scale on the cost function is small and all optimal control strategies, $u(t)$, corresponding to cost functions $c_1(u)$, $c_2(u)$, $c_3(u)$, and $c_4(u)$ appear to suppress infections more aggressively, as illustrated in Figure \ref{fig:lambda1e-2}. The actual values of $\mathcal{I}(t)$ are shown in Table \ref{tab:lambda0.01}. The  total numbers of infected people up to day $t$, $\mathcal{S}(0)-\mathcal{S}(t)$, for all four control functions, $c_1(u)$, $c_2(u)$, $c_3(u)$, and $c_4(u)$, are presented in Table 4.

 \begin{table}[h]
 \caption{Number of infected people $\mathcal{I}(t)$ versus day $t$ for $\lambda=0.01$}
    \centering
    \begin{tabular}{rrrrrr}
    \toprule
       Day & No control & $c_1$ & $c_2$ & $c_3$ & $c_4$ \\
       \midrule
0 & 200 & 200 & 200 & 200 & 200 \\
10 & 1626 & 335 & 340 & 333 & 339 \\
20 & 12542 & 559 & 575 & 550 & 567 \\
30 & 92164 & 932 & 973 & 912 & 951 \\
40 & 643531 & 1563 & 1660 & 1522 & 1599 \\
50 & 2358721 & 2616 & 2829 & 2542 & 2649 \\
60 & 2852657 & 4380 & 4818 & 4245 & 4431 \\
70 & 1722500 & 7447 & 8269 & 7155 & 7587 \\
80 & 834814 & 12778 & 14502 & 12221 & 13336 \\
90 & 373748 & 22520 & 26401 & 21821 & 25197 \\
100 & 164980 & 43921 & 50317 & 42319 & 52546 \\
110 & 71628 & 95049 & 97855 & 87251 & 117329 \\
\bottomrule
    \end{tabular}
    \label{tab:lambda0.01}
\end{table}

\begin{figure}[ht]
    \centering
    \includegraphics[width=.44\textwidth]{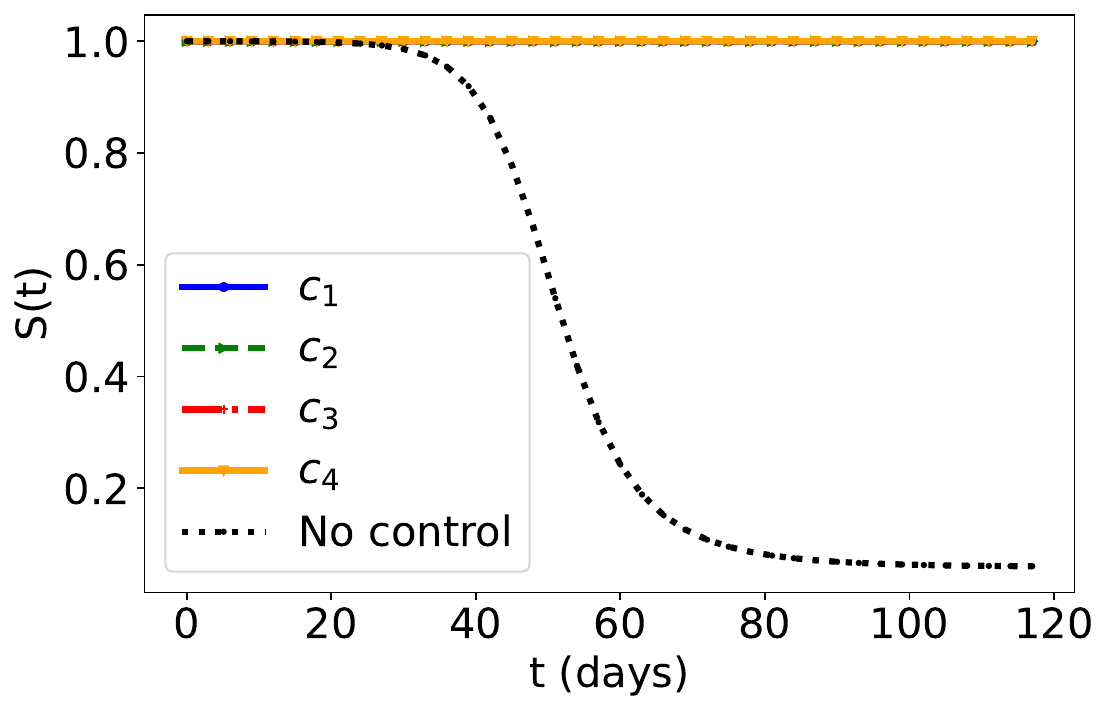}
    \includegraphics[width=.44\textwidth]{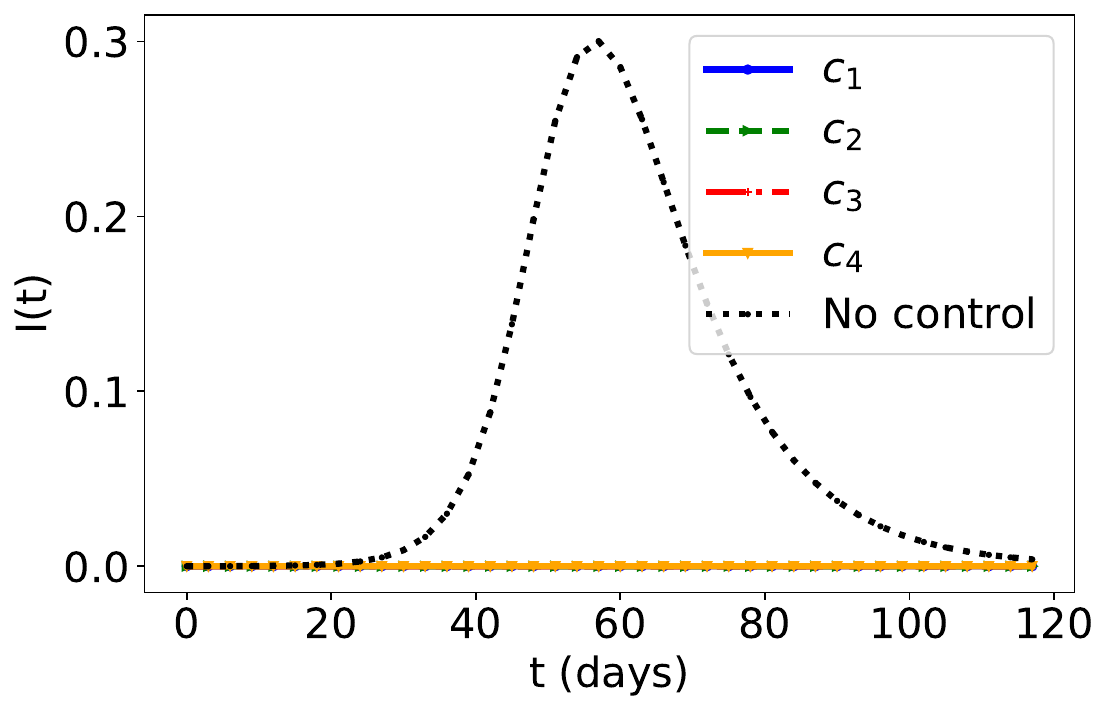}
    \includegraphics[width=.44\textwidth]{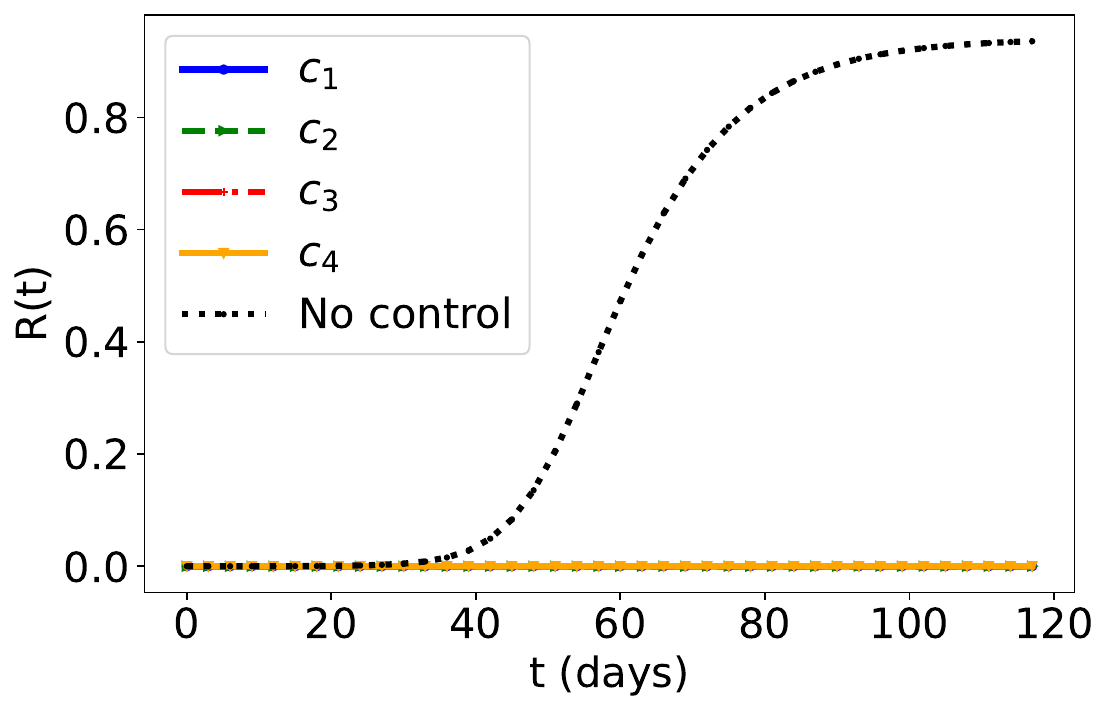}
    \includegraphics[width=.44\textwidth]{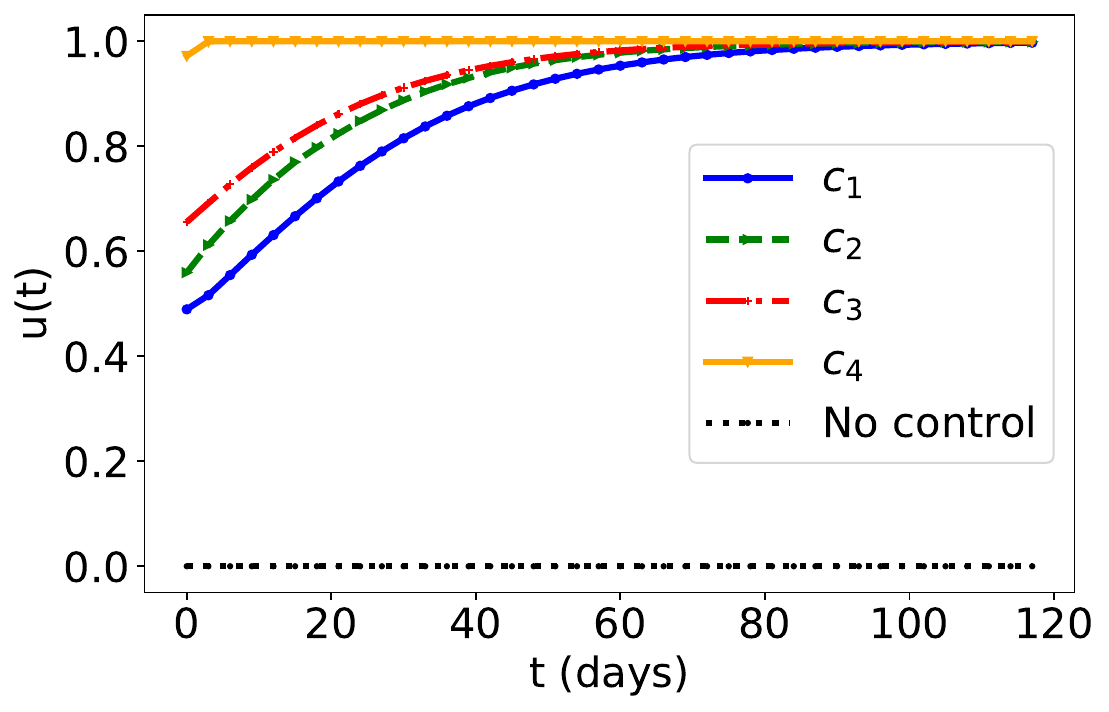}
    \caption{Weight $\lambda=10^{-7}$.}
    \label{fig:lambda1e-7}
\end{figure}

As $\lambda$ continues to decrease, we  see similar behavior as in Figure \ref{fig:lambda1e-2} except that $u(t)$ become larger and the infections are further suppressed. For example, when $\lambda = 10^{-7}$, the cost is very lightly weighted and hence one can impose greater control as shown in Figure \ref{fig:lambda1e-7}.

Again, the behavior of $c_4(u)$ is different. As mentioned above,  this control requires an explicit constraint, $u(t) \in [0,1]$. If not, it is easy to get $u(t)>1$ (especially for small $\lambda$), which is not realistic because it makes $S'(t)>0$. With the constraint enforced, $u(t)$, corresponding to $c_4(u)$, is likely slipping into a local minimum.

As illustrated in  Figure 5, by the time $u(t)$ begins to decrease, the epidemic is effectively under control. Hence, as it follows from (\ref{est}),  the decline in $u(t)$ for $t\ge \tau$ is negligible and the daily number of infected people, $\mathcal{I}(t)$, remains very low for $t\ge \tau$.


For all numerical experiments presented in this section, we let the entire population, $N$  be $10^{7}$ people with $\mathcal{I}(0) = 200$, $\mathcal{R}(0) = 0$, and $\mathcal{S}(0) = N - 200$.

\begin{table}[h]
 \caption{ Total number of infected people up to day $t$,  $N-\mathcal{S}(t)$, for $\lambda=0.01$}
    \centering
    \begin{tabular}{rrrrrr}
    \toprule
Day & No control & $c_1$ & $c_2$ & $c_3$ & $c_4$ \\
\midrule
0 & 200 & 200 & 200 & 200 & 200 \\
10 & 2339 & 604 & 612 & 600 & 615 \\
20 & 18723 & 1276 & 1304 & 1261 & 1289 \\
30 & 138618 & 2394 & 2474 & 2358 & 2437 \\
40 & 993903 & 4287 & 4485 & 4189 & 4369 \\
50 & 4179153 & 7424 & 7899 & 7250 & 7554 \\
60 & 7574716 & 12685 & 13704 & 12335 & 12865 \\
70 & 8800810 & 21705 & 23742 & 20996 & 22079 \\
80 & 9181982 & 37149 & 41310 & 35740 & 38214 \\
90 & 9316270 & 64362 & 73392 & 62132 & 68915 \\
100 & 9367572 & 118524 & 135484 & 114391 & 134473 \\
110 & 9388987 & 236989 & 256274 & 222679 & 281497 \\
\bottomrule
    \end{tabular}
    \label{tab:lambda0.01_N-S}
\end{table}

We use $T=120$, which mimics a 4-months time frame. This value of $T$ allows us to realistically assume that individuals recovered from COVID-19 still have immunity and stay in the removed class, $\mathcal{R}$, for the entire duration of the study period. Furthermore, we take $\beta=0.3$ and $\gamma = 0.1$ days$^{-1}$, which correspond to  the reproduction number, $\mathfrak{R}=3$, and the recovery rate of 10 days.

\section{Conclusions}

In our study, we combine theoretical analysis with rigorous numerical exploration of an optimal control problem for the early stage of an infectious disease outbreak.
We design an objective functional aimed at minimizing the cumulative number of cases. The running cost of control, satisfying  mild, problem-specific, conditions generates an optimal control tragectory, $u(t)$, that stays  inside its admissible set for any $t\in[0,T]$.  For the optimal control problem, restricted by SIR compartmental model (Susceptible $\rightarrow$ Infectious $\rightarrow$ Removed) of disease transmission, we show that the optimal control strategy, $u(t)$,   may be growing until some moment $\tau \in [0 ,T)$. However,  for any $t \in [\tau, T]$,  the derivative, $\dfrac{du}{dt}$, becomes negative and $u(t)$ declines as $t$ approaches $T$ possibly causing the number of newly infected people to go up. So, the window from $0$ to $\tau$ is the time for public health officials  to prepare alternative mitigation measures, such as vaccines, testing,  and antiviral medications,  and to plan for the deployment of rescue equipments like  ventilators and beds.

The impact of $u(t)$ decreasing towards the end of the early stage depends on the weight, $\lambda$. If $\lambda$  is relatively large, then the decline in $u(t)$ for $t\ge \bar{t}$ may be significant, which can result in a considerable surge in the daily number of infected individuals, $\mathcal{I}(t)$, for $t\ge \tau$. On the other hand, if $\lambda$  is small, then  by the time $t=\tau$ the epidemic is effectively under control. Hence, as it follows from (\ref{est}),  the decline in $u(t)$ for $t\ge \tau$ is negligible and the daily number of infected people, $\mathcal{I}(t)$, remains very low for $t\ge \tau$.

Our theoretical findings are illustrated with important numerical examples showing optimal control strategies for various cost functions and weights. Simulation results provide a comprehensive demonstration of the effects of control on the epidemic spread and mitigation expenses, which can serve as invaluable references for public health officials. The important next step is to consider the case of vector-valued controls that, on top of early non-medical interventions (such as social distancing, restrictions on travel and mass gatherings, isolation and quarantine of confirmed cases), include treatment with antivirals and the optimal vaccination strategy.

\section*{Acknowledgement}
We  would like to thank Nathan Gaby for providing assistance to the initial implementation of the code. We would also like to express our deepest gratitude to the associate editor and to the reviewers for their most valuable comments on our manuscript.

\newpage
\bibliographystyle{plain}

\bibliography{refs}

\end{document}